\documentclass[12pt]{article}
\usepackage[left=1in,right=1in,top=1in,bottom=1in]{geometry}

\usepackage{subfigure}
\usepackage{graphicx}
\usepackage{amsfonts}
\usepackage{url}
\usepackage{amsmath,amssymb}
\usepackage{enumerate}
\usepackage{color}
\usepackage{amsthm,amsmath}
\usepackage{bbm}
\usepackage[version=3]{mhchem}
\usepackage{color}

\usepackage{mathtools}

\newcommand{\LRARR}[4]{{\mbox{ \raise 0.4 mm \hbox{$#1$}}} \;
  \mathop{\stackrel{\displaystyle\longrightarrow}\longleftarrow}^{#3}_{#4}
  \; {\mbox{\raise 0.4 mm\hbox{$#2$}}}}

\newcommand{\ra}[1]{\xrightarrow{\hspace*{#1cm}}}
\newcommand{\la}[1]{\xleftarrow{\hspace*{#1cm}}}
\newcommand{\LR}[5]{{\mbox{ \raise 1.3 mm \hbox{$#1$}}} \;
  \mathop{\stackrel{\displaystyle\ra{#5}}{\la{#5}}}^{#3}_{#4} \;
  {\mbox{\raise 1.3 mm\hbox{$#2$}}}}

\newtheorem{definition}{Definition}

\newtheorem{theorem}{Theorem}

\newtheorem{assumption}{Assumption}

\theoremstyle{definition}
\newtheorem{example}{Example}

\newcommand{\R}{{\mathbb R}}
\newcommand{\Z}{{\mathbb Z}}

\newcommand{\E}{\mathbb E}

\newcommand{\Sp}{\mathcal{S}}
\newcommand{\C}{\mathcal{C}}
\newcommand{\Reac}{\mathcal{R}}

\newcommand{\simon}[1]{#1}

\title{Product-form stationary distributions for deficiency zero  networks with non-mass action kinetics}

\author{
David F. Anderson\thanks{Department of Mathematics, University of
  Wisconsin, Madison, USA.  anderson@math.wisc.edu, grant support from
  NSF-DMS-1318832 and Army Research Office grant W911NF-14-1-0401.} \thanks{DFA and SLC would like to thank the Isaac Newton Institute for Mathematical Sciences, Cambridge, for support and hospitality during the program Stochastic Dynamical Systems in Biology: Numerical Methods and Applications, where work on this paper was undertaken. This work was supported by EPSRC grant no EP/K032208/1.},
\and
Simon Cotter\thanks{Department of Mathematics, University of
  Manchester, Manchester, UK. simon.cotter@manchester.ac.uk, grant
  support from EPSRC first grant EP/L023393/1.} \addtocounter{footnote}{-2}\footnotemark}

\begin{document}

\maketitle

\begin{abstract}
  In many applications, for example when computing  statistics of fast subsystems in
  a multiscale setting, we wish to find the stationary
  distributions of systems of continuous time Markov chains. 
    Here we present a class of models that appears naturally in certain averaging approaches whose stationary distributions can be computed explicitly. In particular, we study  continuous time Markov chain models for biochemical interaction systems with non-mass action kinetics whose network satisfies a certain constraint.   Analogous with previous related results, the distributions  can be  written in product form.

%   In this paper, we generalize the results
%  of previous work, to include the form of stationary
%  distributions for systems with the type of non-mass action kinetics
%  that appear when applying certain averaging approaches. 
\end{abstract}


\section{Introduction}
\label{sec:intro}

Biological interaction systems are typically modeled in one of
three ways.  If the counts of the constituent species are
high, then their concentrations are often modeled via a system of
ordinary differential equations with state space $\R^d_{\ge 0}$, where
$d>0$ is the number of species. If the counts are moderate (perhaps order $10^2$ or $10^3$),
  then they may be approximated by some form of continuous diffusion process \cite{gillespie2000chemical,VanKampen}. If the counts are low, then the system is typically modeled stochastically as a continuous time Markov chain in $\Z^d_{\ge 0}$ \cite{AndKurtz2011,AK2015}.   We often want to understand the stationary behavior of the model under consideration.  For deterministic models, understanding the stationary behavior usually entails characterizing the stable fixed points of the system, whereas for stochastic models we require the calculation of the stationary distribution.

  Stationary distributions are also useful in a multiscale setting, where the stationary statistics of a fast subsystem can be utilized in the approximation of the dynamics of the slow variables, which are typically of most interest.  When an analytical form for the stationary distribution of the fast subsystem is not known, numerical approximations can be used.  However, these computations are often expensive and part of an ``inner loop,'' typically making this calculation the rate limiting step of the analysis.   
  %In fact, this  work was motivated by the need to compute invariant distributions of fast subsystems of biochemical processes. 
In the  context of biochemical reaction networks, quasi-equilibrium (QE)-based approximations lead to fast subsystems which
preserve mass action kinetics \cite{goutsias2005quasiequilibrium,janssen1989elimination,thomas2012slow,Cao2005b,weinan2005nested}.  However, 
more recent improvements in stochastic averaging can lead to fast subsystems with non-mass
action
kinetics, and this observation was the motivation for the present work \cite{cotter2015constrained,cotter2014error,cotter2011constrained}.

 One class of interaction networks that has been quite successfully
 analyzed, and that  appears ubiquitously as fast subsystems, are those that are weakly reversible and have a deficiency
 of zero (see Appendix \ref{app:crnt}).  For this class of models, and under the assumption of mass
 action kinetics,   the fixed points of the deterministic models
 \cite{Anderson2011,AndGlobal,Craciun2015,FeinbergLec79,Feinberg87,Gun2003}  and the stationary
 distributions for the stochastic models, \simon{following initial
   work in the 70s~\cite{van1976equilibrium}, have been
 fully characterized~\cite{ACGW2015,ACK2010,CW2016}.}  In fact, it is the study of this class of
 networks that is largely responsible for the development of the field
 of chemical reaction network theory \cite{FeinbergLec79, Gun2003}, a branch of
 applied mathematics in which the dynamical properties of the
 mathematical model are related to the structural properties of the
 interaction network.

In this article, we return to stochastically modeled interaction
networks that are weakly reversible and have a deficiency of zero,
though we consider propensity functions (also called intensity
functions or rate functions) that are more general than
mass action kinetics.  However, we add a \simon{certain condition to
  the rates within the reaction network} (see Assumption \ref{assump:main} below).  Following \cite{ACK2010}, which was motivated
by the work of Frank Kelly \cite{Kelly1979} \simon{who discovered the
  product form stationary distribution of certain stochastically modeled queuing networks}, we provide the form of
the stationary distribution for this  class of models.  In
particular, and in similarity with the main results of \cite{ACK2010},
in Theorem \ref{theoremMain} we show that the distribution is of
product form, and that the key parameter of the distribution is a
complex-balanced equilibrium value of an associated deterministically
modeled system with mass action kinetics. \simon{The result of
  \cite{ACK2010} can now be viewed as special case of this new result,
  where the product form stationary
  distribution is given only in terms of the propensity functions
  evaluated at points in the state space.}

The paper proceeds as follows.  In Section \ref{sec:model}, we introduce the formal mathematical model of interest.    In Section \ref{sec:mainresult} we provide the main theorem of this article, which characterizes the stationary distribution for the class of models of interest.  In Section \ref{sec:examples}, we provide a series of examples which demonstrate the usefulness of the main result. Some brief concluding remarks are given in Section \ref{sec:conc}.   

We assume throughout that the reader is familiar with terminology from chemical reaction network theory.  However, we provide in Appendix \ref{app:crnt}   all necessary terminology and results from this field that are used in the present work.

\section{Mathematical model}
\label{sec:model}

We consider a system with $d$ chemical species, $\{S_1,\dots,S_d\}$,
undergoing reactions which alter the state of the system.  For concreteness, we suppose there are $K>0$ distinct reaction channels.   For the $k$th
reaction channel, we denote by $\nu_k, \nu_k' \in \Z^d_{\ge 0}$ the vectors
representing the number of molecules of each species consumed and
created in one instance of that reaction, respectively.  Note that $\nu_k'-\nu_k\in \Z^d$ is  the net change in the system due to one instance of the $k$th reaction.   We associate each such $\nu_k$ and
$\nu_k'$ with a linear combination of the species in which the
coefficient of $S_i$ is $\nu_{ki}$, the $i$th element of $\nu_k$.  For
example, if $\nu_k = [1, \ 2]^T$ for a system consisting of two
species, then we associate with $\nu_k$ the linear combination $S_1 + 2S_2$.   Under this association, each $\nu_k$ and $\nu_k'$ is
termed a {\em complex} of the system.  We denote any reaction by the
notation $\nu_k \to \nu_k'$, where $\nu_k$ is the source, or reactant,
complex and $\nu_k'$ is the product complex.  We note that each
complex may appear as both a source complex and a product complex in
the system.  The set of all complexes will be denoted by $\{\nu_k\}$.

\begin{definition}
  Let $\Sp = \{S_i\}$, $\C = \{\nu_k\},$ and $\Reac = \{\nu_k \to
  \nu_k'\}$ denote the sets of species, complexes, and reactions,
  respectively.  The triple $\{\Sp, \C, \Reac\}$ is called a {\em
    chemical reaction network}.
  \label{def:crn}
\end{definition}

Throughout we assume that $\nu_k\ne \nu_k'$ for each $k\in \{1,\dots,K\}$.

\subsection{Deterministic model}
The usual deterministic model for a chemical reaction network $\{\Sp,\C,\Reac\}$ assumes that the vector of concentrations for the species satisfies a differential equation of the form
\begin{equation}\label{eq:ode67876}
	\dot x(t) = \sum_{k=1}^K r_k(x(t)) (\nu_k' - \nu_k),
\end{equation}
where $r_k:\R^d_{\ge 0} \to \R_{\ge 0}$ is the (state dependent) rate of the $k$th reaction channel.  If we assume that each function $r_k$ satisfies deterministic mass action kinetics, then
\[
	r_k(x) = \kappa_k \prod_{i = 1}^d x_i^{\nu_{ki}}.
\]

\begin{definition}
	An equilibrium value $c \in \R^d_{\ge 0}$  of \eqref{eq:ode67876} is said to be \emph{complex-balanced} if for each $\eta \in \C$,
	\[
		\sum_{k: \nu_k=\eta} r_k(c) = \sum_{k:\nu_k' = \eta}r_k(c),
	\]
	where the sum on the left, respectively right, is over those reactions with $\eta$ as source complex, respectively product complex.  In the special case of mass action kinetics, $c$ is complex balanced if and only if
	\begin{align}\label{eq:56789875}
	\sum_{k: \nu_k=\eta}\kappa_k \prod_{i = 1}^d c_i^{\nu_{ki}} = \sum_{k:\nu_k' = \eta} \kappa_k \prod_{i = 1}^d c_i^{\nu_{ki}}. 
	\end{align}
\end{definition}

\subsection{Stochastic model, previous results, and assumptions}

The usual stochastic model for a  reaction network $\{\Sp, \C,
\Reac\}$ treats the system as a continuous time Markov chain for which
the rate of transition from state $x\in \Z^d_{\ge 0}$ to state $x +\nu_k'-\nu_k$ is $\lambda_k(x)$, where $\lambda_k:\Z^d_{\ge 0} \to \R_{\ge 0}$ is a suitably chosen intensity function (the intensity functions are also termed \textit{propensity} functions  in the literature).  
%That is, if $X(t)$ is the vector giving the counts of the species at time $t$,
%\begin{align*}
%	P(X(t+\Delta t) = x+\nu_k'-\nu_k \ | \ X(t)=x) &= \lambda_k(x) \Delta t + o(\Delta t)\\
%	P(X(t+\Delta t) = x+\nu_k'-\nu_k \ | \ X(t)=x) &= \lambda_k(x) \Delta t + o(\Delta t)\\
%	P(X(t+\Delta t) = x+\nu_k'-\nu_k \ | \ X(t)=x) &= \lambda_k(x) \Delta t + o(\Delta t),
%\end{align*}
This stochastic process can be characterized  in a variety of useful ways \cite{AndKurtz2011,AK2015}.  For example, it is the stochastic process with state space $\Z^d_{\ge 0}$ and infinitesimal generator 
\begin{equation*}
	Af(x) = \sum_{k=1}^K \lambda_k(x)(f(x+\nu_k'-\nu_k) - f(x)),
\end{equation*}
where $f:\Z^d_{\ge 0} \to \R$.    Kolmogorov's forward equation for this model, termed the chemical master equation in the biology literature, is
\begin{eqnarray}\label{eq:56765}
	\frac{d}{dt} p_\mu(t,x) & = & \sum_{k=1}^K
  \lambda_k(x-\nu_k'+\nu_k) p_\mu(t,x-\nu_k' +
  \nu_k)\mathbbm{1}_{\{x-\nu_k' + \nu_k \ge 0\}} - \sum_{k =1}^K
  \lambda_k(x) p_\mu(t,x), \\
&=& A^*p_\mu(t,x)\mathbbm{1}_{\{x\ge 0\}},
\end{eqnarray}
\noindent where $p_\mu(t,x)$ is the probability the process is in state $x\in \Z^d_{\ge 0}$ at time $t\ge 0$, given an initial distribution of $\mu$, and $A^*$ is the adjoint of $A$.  Note that \eqref{eq:56765} implies that a stationary distribution for the model, $\pi$, must satisfy
\begin{align}
\label{eq:6867767}
 \sum_{k=1}^K \lambda_k(x-\nu_k'+\nu_k) \pi(x-\nu_k' + \nu_k) =\sum_{k =1}^K \lambda_k(x)\pi(x).
\end{align} 
Note also that \eqref{eq:6867767} implies that $\pi$ is in the null space of $A^*$. Thus,  and as is well known, finding $\pi$ reduces to solving $\pi A = 0$, with $\sum_x \pi(x)=1$, where $A$ is reinterpreted as a generator matrix.

%\simon{RANDOM TIME CHANGE PARA REMOVED}

% Another common characterization for this class of models is via the random time change equations of Thomas Kurtz.  Let $N_k(t)$ denote the number of times the $k$th reaction has occurred up to time $t\ge 0$ and let $X(t)$ be the vector whose $i$th coordinate gives the count of species $S_i$ at time $t$.  Then,
% \begin{align*}
% 	X(t) = X(0) + \sum_{k=1}^K N_k(t) (\nu_k' - \nu_k),
% \end{align*}
% and we may represent the counting processes $N_k$ via $N_k(t) = Y_k \left( \int_0^t \lambda_k(X(s)) ds \right)$, where  the $Y_k$ are independent unit Poisson processes \cite{AndKurtz2011,AK2015, Kurtz80}.  This yields  the stochastic equations
% \begin{align*}
% 	X(t) = X(0) + \sum_{k=1}^K Y_k \left( \int_0^t \lambda_k(X(s)) ds \right) (\nu_k' - \nu_k),
% \end{align*}
% which determines the realization $X(t)$ via a jump by jump argument \cite{Anderson2007a}.

Of particular interest to us in this work are the rate functions, $\lambda_k$.   Under the assumption of stochastic mass action kinetics, we have
\begin{equation*}
	\lambda_k(x) =\kappa_k \prod_{i=1}^d \frac{x_i!}{(x_i - \nu_{ki})!} = \kappa_k \prod_{i = 1}^d \prod_{j=0}^{\nu_{ki}-1} (x_i-j),
\end{equation*}
where $\kappa_k>0$ are the \textit{reaction rate constants}.
Here and throughout, we interpret any product of the form $\prod_{i =
  0}^{-1} a_i$ to be equal to one.   In \cite{ACK2010} the authors
considered a class of stochastically modeled reaction networks with
mass action kinetics (those with a deficiency of zero and are weakly reversible, see Appendix \ref{app:crnt}) and characterized their stationary distributions as products
of Poisson distributions.  However, they did  not just consider mass action kinetics in \cite{ACK2010}, but  any kinetics satisfying the functional form
\begin{equation}\label{eq:887555}
\lambda_k(x) = \kappa_k \prod_{i = 1}^d \prod_{j=0}^{\nu_{ki}-1}\theta_i(x_i-j), 
\end{equation}
so long as $\theta_i(0)=0$ for each $i$.  In particular, they proved the following.

\begin{theorem}[Anderson, Craciun, Kurtz, 2010 \cite{ACK2010}]
\label{thm:old}
Let $\{\Sp,\C,\Reac\}$ be a reaction network and for $k\in \{1,\dots,K\}$ let $\kappa_k>0$ be a choice of parameters. Suppose that modeled  deterministically with mass action kinetics and rate constants $\kappa_k$ the system is complex balanced with complex balanced equilibrium $c \in \R^d_{>0}$.  Then, the stochastically modeled system with intensity functions \eqref{eq:887555} admits the invariant measure
	\begin{equation}
	\label{eq:88889}
		\pi(x) = \prod_{i = 1}^d \frac{c_i^{x_i}}{ \prod_{j=0}^{ x_i-1} \theta_i(x_i-j)}.
	\end{equation}
	The measure $\pi$ can be normalized to a stationary distribution so long as it is summable.
\end{theorem}
The connection between weakly reversible, deficiency zero networks and complex balanced equilibria is given in Appendix \ref{app:crnt}.

In this paper, we generalize Theorem \ref{thm:old} by showing how to
find the stationary distributions for models with rates that do not
seem to satisfy the form \eqref{eq:887555}.  However, we add an
assumption pertaining to the \simon{form of the rates within the reaction network,} which
we describe now.  \simon{We begin by partitioning the set of species into two
sets, $\mathcal{S} = \mathcal{S}_1\cup \mathcal{S}_2$.  We will say
that $S_i \in \mathcal{S}_2$ if $\alpha_i = \gcd \{\nu_{1i}, \nu_{2i},
\ldots, \nu_{Ki} \} >1$. Otherwise we say that $S_i \in \mathcal{S}_1$,
noting that $\alpha_i = 1$.}  We now assume that the intensity functions for the reaction network  are  of the form
\begin{align}
\label{eq:rates}
	\lambda_k(x) = \kappa_k \prod_{i = 1}^d \prod_{j=0}^{\frac{\nu_{ki}}{\alpha_i} - 1} \theta_i(x_i - j\alpha_i),
\end{align}
where $\kappa_k>0$ and $\theta_i(x_i) = 0$ if and only if $x_i\le \alpha_i-1$.   

\begin{assumption}
\label{assump:main}
A stochastically modeled reaction network satisfies this assumption if it satisfies the partition described above, and has intensity functions of the form \eqref{eq:rates}.
\end{assumption}

We provide an example to clarify the notation.

\begin{example}
Consider the stochastically modeled system with reaction network
\begin{align}\label{eq:999876}
\ce{$2S_1$ <=>[\lambda_1(x)][\lambda_2(x)] $S_2$}, \quad \ce{$4S_1+2S_2$ <=>[\lambda_3(x)][\lambda_4(x)] $S_3$},
\end{align}
where the intensity functions are placed next to the reaction arrows.  Here $S_1\in \mathcal{S}_2$, with $\alpha_1 = 2$, and $S_2,S_3 \in \mathcal{S}_1$.  The assumption \eqref{eq:rates} then supposes that for appropriate functions $\theta_1,\theta_2,\theta_3:\Z_{\ge 0} \to \R_{\ge 0}$, we have
\begin{align*}
	\lambda_1(x) &= \kappa_1 \theta_1(x_1)\\
	\lambda_2(x) &= \kappa_2 \theta_2(x_2)\\
	\lambda_3(x) &= \kappa_3\theta_1(x_1)\theta_1(x_1-2)   \theta_2(x_2)\theta_2(x_2-1)\\
	\lambda_4(x) &= \kappa_4 \theta_3(x_3).
\end{align*}
  For example, valid choices include $\theta_1(x_1) = x_1(x_1-1)+1(x_1\ge 2)$, $\theta_2(x_2) = x_2$, and $\theta_3(x_3) = \frac{x_3}{1+x_3}$, in which case the form for the stationary distribution does not follow immediately from Theorem \ref{thm:old}.
\end{example}

\section{Main result}
\label{sec:mainresult}

Here we state and prove our main result.  

\begin{theorem}\label{theoremMain}
	Let $\{\mathcal{S}_1\cup\mathcal{S}_2,\C,\Reac\}$ be a reaction network satisfying Assumption \ref{assump:main} and let $(\kappa_1,\dots,\kappa_K)$ be a choice of positive rate constants. Suppose that modeled  deterministically with mass action kinetics and rate constants $\kappa_k$ the system is complex balanced with complex balanced equilibrium $c \in \R^d_{>0}$.  Then the stochastically modeled system with intensity functions \eqref{eq:rates} admits the invariant measure
	\begin{equation}
	\label{eq:88888}
		\pi(x) = \prod_{i = 1}^d \frac{c_i^{x_i}}{ \prod_{j=0}^{\lfloor x_i/\alpha_i\rfloor-1} \theta_i(x_i-j\alpha_i)}.
	\end{equation}
	The measure $\pi$ can be normalized to a stationary distribution so long as it is summable.
\end{theorem}

Note that the theorem applies to models with reaction networks satisfying Assumption \ref{assump:main} and that are weakly reversible and have a deficiency of zero.  See Appendix \ref{app:crnt}. 

\begin{proof}
	First note that if $\mathcal{S}_2 = \emptyset$, then Theorem \ref{theoremMain} is the same as Theorem \ref{thm:old} and there is nothing to show. Thus, we suppose $\mathcal{S}_2 \ne \emptyset$.
	
	The proof proceeds in the following manner.  First, for each $S_i \in \mathcal{S}_2$ we will demonstrate the existence of  a function $\varphi_i:\Z_{\ge 0} \to \R_{\ge 0}$ for which 
	\begin{align}\label{eq:666678}
		\theta_i(x_i) = \prod_{\ell=0}^{\alpha_i-1} \varphi_i(x_i-\ell) =  \varphi_i(x_i)\cdots \varphi_i(x_i-\alpha_i+1).
	\end{align}
	Next, we will apply Theorem \ref{thm:old} and prove that the resulting distribution  is indeed given by \eqref{eq:88888}.
	
	Let $S_i \in \mathcal{S}_2$.  We begin by setting
	\begin{align*}
		\varphi_i(0)=0, \quad \text{ and }\quad  \varphi_i(z) = 1, \text{ for } 1\le z \le \alpha_i-1.
	\end{align*}
	For $z \ge \alpha_i$ an integer,  we may  define $\varphi_i$ recursively via the formula
	\begin{align}\label{eq:3456543}
		\varphi_i(z) = \frac{\theta_i(z)}{\varphi_i(z-1)\cdots \varphi_i(z-\alpha_i+1)}.
	\end{align}
	Note that $\varphi_i$ is a well defined function since $\theta_i(z) > 0$ for each $z \ge \alpha_i$ by assumption.  It is clear that  \eqref{eq:666678} is satisfied with this choice of $\varphi_i$.
	
	For $x\in \Z^d_{\ge 0}$ we may now write
	\begin{align*}
	\lambda_k(x) &= \kappa_k \prod_{i = 1}^d \prod_{j=0}^{\frac{\nu_{ki}}{\alpha_i} - 1} \theta_i(x_i - j\alpha_i) = \kappa_k \prod_{i = 1}^d \prod_{j=0}^{\frac{\nu_{ki}}{\alpha_i} - 1}   \prod_{\ell=0}^{\alpha_i-1} \varphi_i(x_i-j \alpha_i - \ell)= \kappa_k \prod_{i=1}^d \prod_{b=0}^{\nu_{ki}-1} \varphi_i(x_i - b).
	\end{align*}
	Hence, we may apply Theorem \ref{thm:old} and conclude that
	\begin{align*}
		\pi(x) = \prod_{i = 1}^d \frac{c_i^{x_i}}{ \prod_{j=0}^{ x_i-1} \varphi_i(x_i-j)}
	\end{align*}
	is an invariant measure for the system, where $c$ is a complex balanced fixed point for the deterministic system.  It remains to show that
	\begin{equation}\label{eq:54676543}
	\prod_{j=0}^{ x_i-1} \varphi_i(x_i-j) = \prod_{j=0}^{\lfloor x_i/\alpha_i\rfloor-1} \theta_i(x_i-j\alpha_i).
	\end{equation}
	
	First note that if $x_i < \alpha_i$, then both sides of \eqref{eq:54676543} are equal to one.  For the time being, assume that $\alpha_i \le x_i< 2\alpha_i$.  Under this assumption, the right hand side of  \eqref{eq:54676543} is 
	\[
		\theta_i(x_i) = \varphi_i(x_i)\cdots \varphi_i(x_i - \alpha_i + 1) = \varphi_i(x_i)\cdots \varphi_i(\alpha_i ), 
	\]
	where in the final equality we used that $\varphi_i(\ell) = 1$ for $1 \le \ell< \alpha_i$ (and that $x_i - \alpha_i < \alpha_i$), and the left hand side  is 
	\[
		\varphi_i(x_i)\cdots \varphi_i(1) = \varphi_i(x_i)\cdots \varphi(\alpha_i),
	\]
	where we again used that $\varphi_i(\ell) = 1$ when $1\le \ell < \alpha_i$.  Hence, \eqref{eq:54676543} is verified when $x_i <2 \alpha_i$. 
	
	We will now prove that \eqref{eq:54676543} holds in general by induction. We suppose that \eqref{eq:54676543} holds for all $z \le x_i$ where $x_i \ge 2\alpha_i-1$, and will show it to hold at $x_i + 1$.   Using  \eqref{eq:3456543}, the left hand side of \eqref{eq:54676543} evaluated at $x_i+1$ is
	\begin{align*}
		\prod_{j=0}^{(x_i+1)-1}\varphi_i(x_i+1-j) &= \prod_{j=1}^{x_i+1} \varphi_i(j)=\varphi_i(x_i+1)\cdots \varphi_i(x_i+1-\alpha_i + 1) \prod_{j= 1}^{x_i+1 - \alpha_i} \varphi_i(j)\\
		&= \theta_i(x_i+1)\prod_{j= 1}^{x_i+1 - \alpha_i} \varphi_i(j),
	\end{align*}
	where the final equality is an application of \eqref{eq:3456543} with $x_i+1$ in place of the variable $z$.  Continuing, we have
	\begin{align}
	\theta_i(x_i+1)&\prod_{j= 1}^{x_i+1 - \alpha_i} \varphi_i(j)=\theta_i(x_i+1) \prod_{j=0}^{x_i+1-\alpha_i - 1} \varphi_i(x_i+1-\alpha_i - j)\tag{by a rearrangement}\\
	&=\theta_i(x_i+1)  \prod_{j=0}^{\lfloor (x_i+1-\alpha_i)/\alpha_i\rfloor-1} \theta_i((x_i+1-\alpha_i)-j\alpha_i)\tag{by the inductive hypothesis}\\
	&= \theta_i(x_i+1)  \prod_{j=0}^{\lfloor (x_i+1)/\alpha_i\rfloor-2} \theta_i(x_i+1-(j+1)\alpha_i)\notag\\
	&= \theta_i(x_i+1)  \prod_{j=1}^{\lfloor (x_i+1)/\alpha_i\rfloor-1} \theta_i(x_i+1-j\alpha_i)\notag\\
	&= \prod_{j=0}^{\lfloor (x_i+1)/\alpha_i\rfloor-1} \theta_i(x_i+1-j\alpha_i),\label{eq:5676543}
	\end{align}
	where the final equalities are straightforward.   Note that \eqref{eq:5676543} is the right hand side of \eqref{eq:54676543} evaluated at $x_i+1$, and so the proof is complete.
\end{proof}

\section{Examples}
\label{sec:examples}

\subsection{Example 1: Motivating Example}
\label{sec:example1}

First we consider a motivating example arising from 
model reduction, through
constrained averaging \cite{cotter2015constrained,cotter2014error,cotter2011constrained}, of the following system:
\begin{equation}\label{eq:dimer}
\ce{$2S_1$ <=>[\kappa_1x_1(x_1-1)][\kappa_2x_2] $S_2$}, \quad \ce{$\emptyset$ ->[\kappa_3] $S_2$}, \quad \ce{$S_1$ ->[\kappa_4x_1] $\emptyset$},
\end{equation}
where the intensity functions are again placed next to the reaction arrows.
 Note that the intensities of all the reactions follow mass action kinetics.  We consider this system in a parameter regime
where the reversible dimerization reactions $\ce{$2S_1$ <=> $S_2$}$ are
occurring more frequently than the production of $S_2$ and the degradation of
$S_1$. Both $S_1$ and $S_2$ are changed by the fast reactions, but the
quantity $S = S_1 + 2S_2$ is invariant with respect to the fast reactions,
and as such is the slow variable in this system. We wish to reduce the
dynamics of this system to a model only concerned with the possible changes in $S$:

\begin{equation}\label{eq:dimer:eff} 
\ce{$\emptyset$ ->[\bar{\lambda}_3(s)] $2S$}, \quad \ce{$S$ ->[\bar{\lambda}_4(s)] $\emptyset$},
%\Ra{\emptyset}{2S}{\bar{\lambda_3}(S)}{2}, \qquad \Ra{S}{\emptyset}{\bar{\lambda_4}(S)}{2}
\end{equation}
where $\bar{\lambda}_3(s)$ and $\bar{\lambda}_4(s)$  are the effective rates of the system.

Using the QE approximation (QEA),  $\bar{\lambda}_3(s)=\kappa_3$ and  %distribution $\mathbb{P}(Y|S)$, which we require in
 $\bar{\lambda}_4(s)=\kappa_4 \E_{\pi_{\rm QEA}(s)}[X_1]$, where $\pi_{\rm QEA}(s)$ is the stationary distribution for the system
%   the mean value  by the
%invariant distribution of the following system:
\begin{equation}\label{eq:dimer:fss}
\ce{$2S_1$ <=>[\kappa_1x_1(x_1-1)][\kappa_4x_2] $S_2$}, 
%\LR{2X}{Y}{k_3X(X-1)}{k_4Y}{3}, \qquad X + 2S_2 = S.
\end{equation}
under the assumption that $X_1(0) + 2X_2(0) = s$.
Since the system \eqref{eq:dimer:fss} satisfies the necessary conditions of the results of
\cite{ACK2010} (weak reversibility and deficiency of zero), the invariant distribution $\pi_{\rm QEA}(s)$ is known exactly.

In comparison, the constrained approach requires us to find the
invariant distribution $\pi_{\rm Con}$ of the following system:
\begin{equation}\label{eq:con}
\ce{$2S_1$ <=>[\kappa_1x_1(x_1-1)+\kappa_3\mathbbm{1}_{\{x_1>1\}}][(\kappa_2+\kappa_4)x_2] $S_2$},
\end{equation}
subject to $X_1(0)+2X_2(0)=s$.
Readers interested in seeing how this is derived should refer to \cite{cotter2015constrained}.
This network is weakly reversible, and has a deficiency of zero. However,
the form of the rates in this system do not satisfy the conditions
specified in \cite{ACK2010}.
%, and as such we do not currently have a closed form for the stationary distribution of this simple system.
\simon{ In
the context of constrained averaging, this lack of a closed form for the stationary distribution would result in the need for some form of approximation of the stationary distribution.  There are two common methods utilized for performing this approximation.  One possibility would be to perform exhaustive stochastic simulation of the system \eqref{eq:con}.  Another option involves finding the distribution by finding the null space of the adjoint of the generator (see the discussion in and around \eqref{eq:6867767}).  However, as the state space of  \eqref{eq:con} will typically be huge, the latter method often involves truncating the state space and approximating the actual distribution with that of the stationary distribution of the truncated system  \cite{cotter2015constrained}.} Both approaches will lead to approximation errors
and varying amounts of computational cost. %This motivates the need
                                %for development of theory regarding
                                %stationary distribution of stochastic
                                %chemical reaction networks with
                                %non-mass action kinetics.
However, note that the system \eqref{eq:con} does satisfy Assumption \ref{assump:main}, with $\alpha_1 = 2$ and $\alpha_2 = 1$. We denote the
rate of dimerization by $\lambda_D$ and its reverse by
$\lambda_{-D}$. Therefore
\begin{eqnarray*}
\lambda_D(x) &=& k_1x_1(x_1-1) + k_3\mathbbm{1}_{\{x_1>1\}},\\
&=& k_1\left (
  x_1(x_1-1) + \frac{k_3}{k_1}\mathbbm{1}_{\{x_1>1\}} \right ),\\
&=& k_1 \theta_1(x_1),
\end{eqnarray*}
with $\theta_1$ defined in the final equality.
The form of the rate of the reverse reaction is much simpler, and is given by
\begin{eqnarray*}
\lambda_{-D}(x) &=& (k_2+k_4)x_2 = (k_2 + k_4)\theta_2(x),
\end{eqnarray*}
which defines $\theta_2$.

By Theorem \ref{theoremMain}, we can write down the stationary
distribution of this system. The complex balanced equilibrium of the associated
deterministically modeled system with $c_1+2c_2=1$, is given by $(c_1, c_2) = \left
  (\frac{ \sqrt{(k_2+k_4)(k_2 + 8k_1 + k_4)} - k_2 - k_4}{4k_1},
  \frac{1 - c_1}{2} \right)$.  Then by Theorem \ref{theoremMain}, and
by recalling that all states $(x_1,x_2)$ in the domain satisfy $s = x_1 + 2x_2$, the
stationary distribution for $S_2$ is given by
\begin{equation}\label{eq:dimer:con}
\pi_{\rm{Con}}(x_2) = \frac{1}{\Gamma_{\rm{Con}}}
\frac{c_1^{(s-2x_2)}}{\prod_{j=0}^{\lfloor (s-2x_2)/2\rfloor-1}
   \left ((s-2x_2-2j)(s-2x_2 -2j-1)
  + \frac{k_3}{k_1} \right )} \frac{c_2^{x_2}}{x_2!},
\end{equation}
where $\Gamma_{\rm{Con}}$ is a normalization constant and $s$ is the conserved quantity. Note that the
indicator function in $\theta_1$ (in the denominator) has disappeared since it is always equal to one over the
domain of the product.

We can compare \eqref{eq:dimer:con} with the distribution of \eqref{eq:dimer:fss}, which arises from the QEA, and also with the distribution of the full system \eqref{eq:dimer} conditioned on $S_1 + 2S_2 =s$ (which can be approximated
  by finding the null space of the adjoint of the  generator of the full system on the truncated domain).
First we consider the QEA approximation. The invariant distribution
of the fast subsystem \eqref{eq:dimer:fss} can be found using Theorem
\ref{thm:old}, and is given by
\begin{equation}\label{eq:dimer:qssa}
\pi_{\rm{QEA}}(x) = \frac{1}{\Gamma_{\rm{QEA}}} \frac{d_1^{s-2x_2}}{(s-2x_2)!} \frac{d_2^{x_2}}{x_2!},
\end{equation}
where $(d_1, d_2) = \left ( \frac{\sqrt{k_4(8k_1 +
      k_4)}-k_4}{4k_1},\frac{1- d_1}{2} \right ) $ is the complex
balanced equilibrium for this system satisfying $d_1+2d_2=1$, 
and $\Gamma_{\rm{QEA}}$ is a normalizing constant.

Since the full system \eqref{eq:dimer} does not have a deficiency of zero,
we are not able to find its invariant distribution directly. However,
by truncating the state space appropriately, we are able to
approximate the full distribution by constructing the generator on
this truncated state space and finding the null space of the adjoint.
% \simon{ The direction of the one-dimensional null space, once normalized, provides us with the stationary distribution of the truncated process.} 

%The generator
%$\mathcal{G}$ is the matrix which allows us to put the chemical master
%equation \eqref{eq:56765}
%into matrix-vector format given by
%\begin{equation*}
%\frac{d {\bf P}}{dt} = \mathcal{G},
%\end{equation*}
%where ${\bf P}$ is the vector of probabilities $\mathbb{P}(x,t)$ of
%being in state $x \in \mathbb{Z}^d_{\ge 0}$ at time $t$, and therefore the invariant
%distribution is given by the null space of the generator.
Once we have approximated the null space of the truncated generator, we
can find the approximation of $\mathbb{P}(X_2 = x_2|X_1 + 2X_2=s)$ by taking the
probabilities of all states with $x_1+2x_2 = s$ and renormalizing. In what
follows, we truncated the domain of the generator to $x \in
\{0,1,\ldots,1000\} \times \{0,1,\ldots,500\}$.

We consider the system \eqref{eq:dimer} with parameters given by:
\begin{equation}\label{eq:dimer:params}
k_1 = 1, \qquad k_2 = 100, \qquad k_3 = 1500, \qquad k_4 = 30.
\end{equation}
Note that it is not obvious from these rates that the reactions with
rates $k_3$ and $k_4$ are in fact the slow reactions in this
system. The invariant density is largely concentrated in a small
region centered close to the point $x = (99,114)$. By using the
approximation of the invariant density that we have computed on the
truncated domain, we can compute the expected ratio between
occurrences of the fast reactions with rates $k_1$ and $k_2$ with the
slow reactions with rates $k_3$ and $k_4$. For this choice of
parameters, the expected proportion of the total reactions which are
fast reactions (dimerization/disassociation) is $82.68\%$. This indicates a difference in
timescales between these reactions, but the difference is not particularly
stark, and as such we would expect there to be significant error in
any approximation relying on the QEA.

\begin{figure}
\begin{center}
\includegraphics[width=0.5\textwidth]{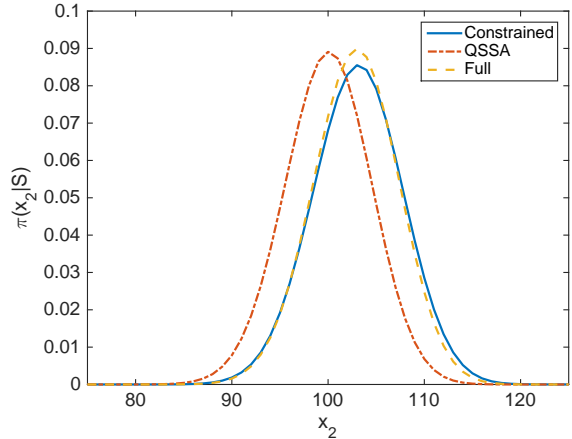}
\caption{Approximations of the distribution $\mathbb{P}(X_2=x_2|X_1 + 2X_2=300)$ for
  the system \eqref{eq:dimer} with parameters given by
  \eqref{eq:dimer:params} using constrained averaging, QEA
  averaging, and through approximation of the invariant distribution
  of the full system on $x \in
\{0,1,\ldots,1000\} \times \{0,1,\ldots,500\}$. \label{Fig:Example1}}
\end{center}
\end{figure}

Figure \ref{Fig:Example1} shows the three approximations of  the
distribution $\mathbb{P}(X_2=x_2|X_1+2X_2=300)$ for  the system \eqref{eq:dimer} with parameters given by
  \eqref{eq:dimer:params}. The constrained and QEA approximations are computed
  using \eqref{eq:dimer:con} and \eqref{eq:dimer:qssa} respectively,
  with the normalizing constants computed numerically. As would be
  expected in this parameter regime, the constrained approximation is
  far more accurate than the QEA. 

We can quantify the accuracy of each of the approximations by
computing the relative $l^2$ differences with the distribution computed
using the full generator on the truncated domain. This relative
difference was $4.464 \times 10^{-1}$ for the QEA, in comparison with
$5.2337\times 10^{-2}$ for the constrained approximation.
% Kullback-Leibler (KL) divergences
%   of the constrained and QSSA distributions ($\pi_{\rm Con}$ and
%   $\pi_{\rm QSSA}$) with respect to the
%   distribution computed using the full generator on the truncated
%   domain. The KL
%   divergence is a way to quantify how different two probability
%   distributions are, and for a discrete distribution such as these, on
%   the domain $0 \leq x_2 \leq \lfloor s/2 \rfloor$, it
%   is given by:
% \begin{equation*}
% D_{\rm KL}(\pi_1||\pi_2)= \sum_{0\le x_2 \le \lfloor s/2 \rfloor}
% \pi_1(x_2) \log \left (\frac{\pi_1(x_2)}{\pi_2(x_2)} \right ).
% \end{equation*}
% Two identical distributions will give a KL divergence of
% 0. The further this quantity is from 0, the bigger the difference
% between the distributions.
% For the example shown in Figure \ref{Fig:Example1}, the constrained
% approximation has a KL divergence of 1.0091, compared with the QSSA
% approximation which has a KL divergence of 1.5135. 
This demonstrates
the improvement in approximation that can be achieved by using
constrained averaging, and which motivates the need for results like Theorem \ref{theoremMain} which take non-mass action kinetics into account.

\subsection{Example 2: Dimerization}
\label{sec:dimer}

We begin by considering a model consisting of only proteins, denoted $P$, and dimers, denoted $D$.  We suppose there are two mechanism by which the proteins are dimerized: random interactions between the protein molecules, and via a catalyst.  The rate at which random interactions lead to the formation of dimers can be taken to be of mass action form.  Assuming the concentrations are such that the catalyst is acting at capacity,  the rate of formation of the dimers due to the catalyst can be faithfully modeled as a constant (so long as there are at least two proteins to make the reaction happen).  Thus, letting $x_p$ and $x_d$ denote the numbers of proteins and dimers in the model, respectively, the reaction system can be represented as
\begin{equation}\label{eq:5676546700}
\ce{$2P$ <=>[\kappa_{p\to d}(x_p(x_p-1)+\rho 1(x_p>1))][\kappa_{d\to{p}}x_d] $D$},
\end{equation}
where $\kappa_{p\to d}$ and $\rho$ are given parameters.  Note that after an obvious change of variables the stationary distribution of this particular model is provided in \eqref{eq:dimer:con}.  

The reactions in \eqref{eq:5676546700} are typically a subset of the reactions in a larger system.   For example, the actual model of interest may be
\begin{align}
\ce{$\emptyset$ <=>[\kappa_m][d_mx_m] $M$ ->[\kappa_{m\to p}x_m] $P$ ->[d_px_p] $\emptyset$},\label{eq:5677544}\\
\ce{$2P$ <=>[\kappa_{p\to d}(x_p(x_p-1)+\rho1(x_p>1))][\kappa_{d\to p}x_d] $D$ <=>[d_Dx_d][\kappa_d] $\emptyset$},\nonumber
\end{align}
where $M$ represents an mRNA molecule.  This is a standard model for dimer production.  Depending upon the  relevant time-scales in the system we may want to take $d_D = \kappa_d = 0$.  If the reactions $2P\rightleftharpoons D$ are appreciably faster than those of \eqref{eq:5677544}, then an obvious path simulation strategy presents itself in the mold of the slow scale SSA \cite{Cao2005b} (i.e.~stochastic averaging):
\begin{enumerate}
\item for the current values of $P$ and $D$, find the stationary distribution of the fast subsystem analytically,
\item determine the effective rates of the reduced model
\begin{align}
\ce{$\emptyset$ <=>[\kappa_m][{d_m}{x_m}] $M$ ->[\kappa_{m\to{p}}x_m] $P$ ->[d_p\E(x_p)] $\emptyset$},  \quad \ce{$D$ <=>[d_DE(x_d)][\kappa_d] $\emptyset$}
\end{align}
where the expectations are with respect to the distribution found in step 1.
\item simulate forward in time using the stochastic simulation algorithm \cite{Gill76} or the next reaction method \cite{Anderson2007a,Gibson2000}, and return to step 1.
\end{enumerate}
Being able to analytically calculate the stationary distribution in step 1 allows us to bypass the need to numerically approximate the stationary distribution, as is commonly done \cite{weinan2005nested,EVE2007}.

\subsection{Example 3}
\simon{In Sections \ref{sec:example1} and \ref{sec:dimer}, we applied Theorem \ref{theoremMain} on the common motif $2S_1 \rightleftharpoons S_2$.  In this example, we present another common motif, $2S \rightleftharpoons \emptyset$,  for which Theorem \ref{theoremMain} is also useful.  As opposed to the specific models we considered in the previous examples, here we present a more general framework in which the specific form of the propensity function for the reaction $2S \to \emptyset$ is arbitrary.  This situation is common when undertaking certain types of
  averaging arguments \cite{cotter2015constrained}.}
  
  \simon{Let us suppose that the effective dynamics of a slow variable in a larger system can be modeled as
%we wish to find the invariant distribution of a
%slow variable $S$ in a system for which we have estimated the
%effective dynamics using averaging. Specifically, we want the invariant distribution of the following model,
\begin{equation*}
\ce{$\emptyset$ ->[k^b] $2S$ ->[k^d_1\mathbb{E}(f_1,\dots,f_r|S=s)] $\emptyset$},
%\Ra{\emptyset}{2S}{k^b}{1.5}\Ra{}{\emptyset}{k^d_1\mathbb{E}(f_1|S)+k^d_2\mathbb{E}(F_2|S)}{3.5},
\end{equation*}
where $k^b$ and $k_1^d$ are positive constants and where $\mathbb{E}(f_1,\dots,f_r|S=s)$ is a conditional expectation of the fast variables $f_1,\dots,f_r$, conditioned on $S=s$.  The conditional expectation could in general be highly nonlinear, and not of the form \eqref{eq:887555} required by Theorem \ref{thm:old}.
Supposing that $\mathbb{E}(f_1,\dots,f_r|S=s) = 0$ if $s\in \{0,1\}$, Theorem \ref{theoremMain} says that, up to a
normalization constant, the invariant distribution of $S$ is given by
\begin{equation*}
\pi(s) \propto \frac{ \left (\frac{k^b}{k^d_1} \right
  )^s}{ \prod_{j=0}^{\lfloor s/2\rfloor-1}  \mathbb{E}(f_1,\dots,f_r|S=s-2j) }.
\end{equation*}
 In practice the normalization
constant could be approximated by summing over an
appropriate domain.
}

\subsection{Example 4}
  This example will demonstrate the difficulties that can arise when a
reaction is added that involves a  species in the set $\mathcal{S}_2$ with multiplicity not equal to $\alpha_i$. 
In particular, consider the system
\begin{equation}\label{eq:2X1X}
\ce{$2S_1$ <=>[\lambda_1(x)][\lambda_2] $\emptyset$}, \quad \ce{$S_1$ <=>[\lambda_3(x)][\lambda_4(x)] $S_2$},
\end{equation}
where
\begin{eqnarray*}
\lambda_1(x) = \kappa_1 \theta_1(x_1), \quad
\lambda_2(x) = \kappa_2 , \quad
\lambda_3(x) = \kappa_3 \varphi(x_1),\quad
\lambda_4(x) = \kappa_4x_2,
\end{eqnarray*}
where $\theta_1(x_1) = \left (\mathbbm{1}_{\{x_1>1\}}(10 + x_1 +
  6\sin(\pi x_1/5)) \right )$.
Note that we have a non-mass action kinetics rate $\lambda_1$ for the
 reaction  $2S_1 \to \emptyset$. There is another reaction
involving $S_1$, but the amount of $S_1$ molecules involved in this
reaction is not a multiple of 2. This means that we cannot apply the
result of Theorem \ref{theoremMain} to this system unless $\varphi$ satisfies the recurrence relation detailed in the proof of Theorem \ref{theoremMain}.  In particular, it must satisfy
$\varphi(x_1)\varphi(x_1-1) = \theta_1(x_1)$, or
\begin{eqnarray*}
\varphi(0) &=& 0, \\
\varphi(1) &=& C \in \mathbb{R}_{>0},\\
\varphi(n) &=& \frac{\theta_1(n)}{\varphi(n-1)}.
\end{eqnarray*}
This recurrence relation defines a unique function
$\varphi:\mathbb{Z}_{\ge 0} \to \mathbb{R}_{\ge 0}$ for each $C \in
\mathbb{R}_{\ge 0}$. 
%Assuming that $\theta_1$ grows asymptotically at a rate which  guarantees the system has a stationary distribution,  then 
In
general the function $\varphi$ can oscillate wildly. Let us consider,
for example, $\varphi$ when $\theta_1$ is given as above.  In this case
\begin{equation}\label{eq:2X1X:varphi}
\varphi(x_1) = \frac{ C^{2{\rm mod}(x_1,2) - 1}\prod_{i=0}^{\lfloor x_1/2 \rfloor -
    1}\theta_1(x_1 - 2i)}{ \prod_{i=0}^{\lfloor (x_1-1)/2 \rfloor -
    1}\theta_1(x_1 - 2i - 1)}.
\end{equation}

\begin{figure}
\begin{center}
\includegraphics[width=0.5\textwidth]{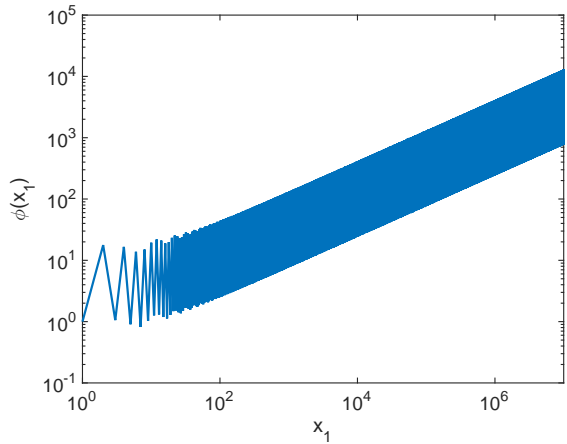}
\caption{$\varphi$ as given in \eqref{eq:2X1X} with $C=1$.\label{Fig:2X1X:varphi}}
\end{center}
\end{figure}
Figure \ref{Fig:2X1X:varphi} demonstrates how $\varphi(x_1)$, and the
amplitude of the oscillations, grow with $x_i$, for the case $C=1$. It is clear that this
function does not represent any physical reaction rate arising from
chemistry, but we can still write down the invariant distribution of
the system \eqref{eq:2X1X}. 

%All species are in $\mathcal{S}_1$ and as
%such for this system the result in Theorem \ref{theoremMain} agrees
%with the result in \cite{ACK2010}.

The complex balance equilibrium of the associated mass action
kinetics system to \eqref{eq:2X1X} is given by $(c_1,c_2) = \left(
  \sqrt{\frac{k_2}{k_1}}, \frac{k_3}{k_4}\sqrt{\frac{k_2}{k_1}} \right
)$. Therefore, the stationary distribution is given by:
\begin{equation}\label{eq:2X1X:SD}
\pi(x) = \frac{1}{\Gamma} \frac{c_1^{x_1}}{\prod_{i=1}^{x_1}
  \varphi(x_1)} \frac{c_2^{x_2}}{x_2!},\qquad x \in \{(x_1,x_2)\ |\  {\rm
  mod}(x_1+x_2,2) = a\} \quad {\rm for} \quad a \in \{0,1\},
\end{equation}
where $\Gamma$ is a normalizing constant, and $\varphi$ is given by
\eqref{eq:2X1X:varphi} with $C=1$.  Note that the value of $a$ here dictates the
oddness or evenness of the quantity $x_1+x_2$, which is preserved by each of the  reactions.

\begin{figure}
\begin{center}
\includegraphics[width=0.5\textwidth]{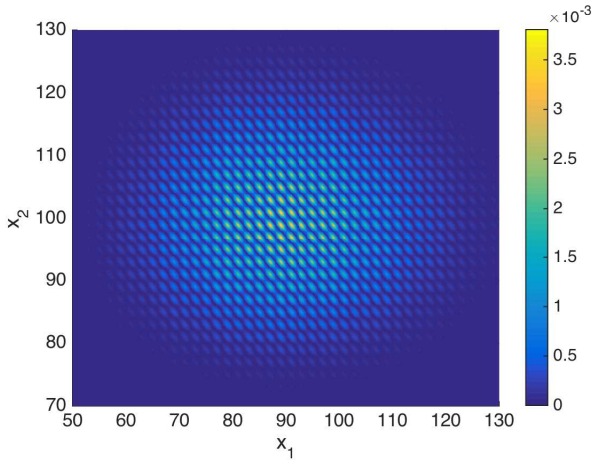}
\caption{Stationary distribution \eqref{eq:2X1X:SD} of the system
  \eqref{eq:2X1X} with parameters given by
  \eqref{eq:2X1X:params}, with the assumption that the initial
  value of $x_1+x_2$ is even. The normalization constant $\Gamma$ was
  approximated by summing the value of all states $x \in
  \{0,1,\ldots,1000\} \times \{0,1,\ldots,1000\}$.  \label{Fig:2X1X}}
\end{center}
\end{figure}

Figure \ref{Fig:2X1X} shows the invariant distribution
\eqref{eq:2X1X:SD} of the system \eqref{eq:2X1X} with an even initial
value of $x_1 + x_2$, and for the following
parameters:
\begin{equation}\label{eq:2X1X:params}
k_1 = 1, \qquad k_2 = 10^2, \qquad k_3 = 10, \qquad k_4 = 1.
\end{equation}
The normalization constant was approximated by summing the values of
all states $x \in
  \{0,1,\ldots,1000\} \times \{0,1,\ldots,1000\}$.   
%  The relative
%  difference between this solution and the numerically found
%  null space of the generator of the process truncated to the same
%  domain was $7.3130\times 10^{-14}$. The equivalent relative error
%  for the solution on the same domain but with an odd initial value of
%  $x_1 + x_2$ was $7.1084\times 10^{-14}$. 

\section{Conclusions}
\label{sec:conc}

 In this paper we provided the stationary distributions for the class of stochastically modeled reaction networks with non-mass action kinetics that satisfy our Assumption \ref{assump:main} and that admit a complex balanced equilibrium when modeled deterministically with mass action kinetics.  Similarly with the results of \cite{ACK2010}, we showed that the stationary distributions are of product form.

We motivated the need for such results through consideration of modern averaging techniques.
%In particular, four examples were presented with two demonstrating how the main result can be used in two different ways when considering averaging of multiscale networks. 
%The last demonstrated the main result's
%application to even the most abstract and unphysical of networks.
In particular, Theorem \ref{theoremMain} significantly  reduces  the computational cost of finding
the invariant distribution of most fast subsystems in a multiscale setting,  therefore making
accurate  approximations very cheap to compute. This in turn
opens up more possibilities, for example approximation of likelihoods
via multiscale reductions in the context of parameter inference for biochemical
networks. 

\appendix

\section{Terminology and results from chemical reaction network theory}
\label{app:crnt}

For each reaction network, $\{\Sp,\C,\Reac\}$, there is a unique directed
graph constructed in the following manner.  The nodes of the graph are
the complexes, $\C$.  A directed edge is then placed from complex
$\nu_k$ to complex $\nu_k'$ if  $\nu_k \to \nu_k' \in \Reac$.
Each connected component of the resulting graph is termed a {\em
  linkage class}.  We denote the number of linkage
classes by $\ell$.  For example, in the reaction network \eqref{eq:999876} there are two linkage classes.

\begin{definition}
  A chemical reaction network, $\{\Sp,\C,\Reac \}$, is called {\em weakly
    reversible} if each  linkage classes is strongly connected.  A network
  is called {\em reversible} if $\nu_k' \to \nu_k \in \Reac$ whenever
  $\nu_k \to \nu_k' \in \Reac$.
  \label{def:WR}
\end{definition}

Note that a network is weakly reversible if and only if for any reaction $\nu_k \to \nu_k'$, there is a
  sequence of directed reactions beginning with $\nu_k'$ as a source
  complex and ending with $\nu_k$ as a product complex.  That is,
  there exist complexes $\nu_1,\dots,\nu_r$ such that $\nu_k' \to
  \nu_1, \nu_1 \to \nu_2, \dots, \nu_r \to \nu_k \in \Reac$.

\begin{definition}
  $S = \hbox{\rm{span}}_{\{\nu_k \to \nu_k' \in \Reac\}}\{\nu_k' - \nu_k\}$
  is the {\em stoichiometric subspace} of the network. For $z \in
  \R^d$ we say $z + S$ and $(z + S) \cap \R^d_{>0}$ are the {\em
    stoichiometric compatibility classes} and {\em positive
    stoichiometric compatibility classes} of the network,
  respectively.  Denote $\hbox{\rm{dim}}(S) = s$.
\end{definition}

%It is clear that for both stochastic and deterministic models, the state of the system remains within a single stoichiometric compatibility class for all time. 

The final definition  is that of the {\em deficiency} of a network \cite{Feinberg87}. 
\begin{definition}
  The {\em deficiency} of a chemical reaction network,
  $\{\Sp,\C,\Reac\}$, is $\delta = |\C| - \ell - s$, where $|\C|$ is the
  number of complexes, $\ell$ is the number of linkage classes of the
  network graph, and $s$ is the dimension of the stoichiometric
  subspace of the network.
\end{definition}

  We state a classical result that can be found in  \cite{FeinbergLec79,Feinberg87, Gun2003}, which  relates networks that are weakly reversible and have a deficiency of zero to those that admit complex balanced equilibria.
%  \begin{theorem}\label{thm:complex_balanced}
%  If a deterministic mass action system  possesses a complex balanced equilibrium, then the system is complex balanced and $\G$ is weakly reversible. Moreover, there exists exactly one complex balanced equilibrium in every positive stoichiometric compatibility class, and it is locally asymptotically stable relative to its positive stoichiometric compatibility class. 
% \end{theorem}
 \begin{theorem}\label{thm:deficiency_zero_iff}
  If the reaction network $\{\Sp,\C,\Reac\}$ is weakly reversible and has a deficiency of zero, then for any choice of rate constants $\kappa_k$ the deterministic mass action  system admits a complex balanced equilibrium, $c$, satisfying \eqref{eq:56789875}.  Moreover, within each positive stoichiometric compatibility class there is precisely one equilibrium, and it is complex balanced.
 \end{theorem}

 \bibliography{bib}
 \bibliographystyle{plain}
 
\end{document}